\documentclass[12pt,headings=normal]{scrartcl}

\usepackage[utf8]{inputenc}
\usepackage[ngerman]{babel}
\usepackage{amsmath,amssymb}
\usepackage[thmmarks]{ntheorem}
\usepackage{verbatim}
\usepackage{graphicx}
\usepackage[all]{xy}

\theorembodyfont{\itshape}
\theoremstyle{plain}
\newtheorem{Satz}{Satz}[section]
\newtheorem{Lemma}[Satz]{Lemma}
\newtheorem{Folgerung}[Satz]{Folgerung}
\newtheorem{Proposition}[Satz]{Proposition}
\theorembodyfont{\upshape}
\newtheorem{Bemerkung}[Satz]{Bemerkung}

\theoremheaderfont{\sc}
\theorembodyfont{\upshape}
\theoremstyle{nonumberplain}
\theoremseparator{}
\theoremsymbol{$\Box$}
\newtheorem{Beweis}{Beweis}

\title{\Large \"Uber die von einem Ideal $I \subset R$ erzeugten $R$-Moduln}
\author{\large Helmut Zöschinger\\
  \large Mathematisches Institut der Universität München\\
  \large Theresienstr. 39, D-80333 München\\
  \large E-mail: zoeschinger$@$mathematik.uni-muenchen.de
}
\date{}

\newcounter{abccount}
\newenvironment{abc}{%
\begin{list}{(\alph{abccount})}{%
\usecounter{abccount}%
\setlength{\partopsep}{0pt}%
\setlength{\topsep}{1ex}%
\setlength{\itemsep}{0pt}%
}%
}{\end{list}}

\newcounter{myenumcount}
\newenvironment{myenum}{%
\begin{list}{(\arabic{myenumcount})}{%
\usecounter{myenumcount}%
\setlength{\partopsep}{0pt}%
\setlength{\topsep}{1ex}%
\setlength{\itemsep}{0pt}%
}%
}{\end{list}}

\newcounter{iiicount}
\newenvironment{iii}{%
\begin{list}{(\roman{iiicount})}{%
\usecounter{iiicount}%
\setlength{\labelwidth}{3em}%
\setlength{\partopsep}{0pt}%
\setlength{\topsep}{1ex}%
\setlength{\itemsep}{0pt}%
}%
}{\end{list}}

\begin{document}
\maketitle

\centerline{\textbf{Abstract}}
\begin{abstract}
  \noindent
  Let $(R, \mathfrak m)$ be a commutative noetherian local ring. We investigate
  under which conditions an $R$-module $M$ is generated by an ideal $I$,
  i.e. there exists an epimorphism $I^{(\Lambda)} \twoheadrightarrow M$. If
  $M$ is uniserial, i.e. $\mathcal{L}(M)$ is totally ordered and finite,
  this is equivalent to $\mathfrak{m}^{n-1} \cdot I \not\subset
  \operatorname{Ann}_R(M) \cdot I$ ($\operatorname{length}(M) = n \geq 1$).
  If $M$ is cyclic and $I = \mathfrak{m}$, this is equivalent to: Either it
  is $M \cong R/\mathfrak{p}$ ($R/\mathfrak{p}$ a discrete valuation ring)
  or $M \cong C/\operatorname{So}(C)$ ($C$ a uniserial $R$-module). If $A$
  is free and $B$ is a submodule of $A$, then the Matlis dual $(A/B)^{\circ}
  = \operatorname{Hom}_R(A/B, E)$ is $I$-generated if and only if $B = (IB)
  :_A I$. In the case $I = \mathfrak{m}$, this condition leads to the
  ``basically full ideals'' considered by Heinzer, Ratliff~Jr. and Rush. By
  studying the dual condition $M = I(M :_X I)$ in the last section, we can
  generalize some results of that work.
\end{abstract}

\bigskip


\bigskip


\section{$I$-generierte Moduln}

Sei $(R, \mathfrak m)$ ein kommutativer, noetherscher lokaler Ring, $I$ ein
Ideal von $R$ und $M$ ein $R$-Modul. Bekanntlich hei"st $M$
\emph{$I$-generiert} (in Zeichen $M \in \operatorname{Gen}(I)$), wenn es
eine Indexmenge $\Lambda$ und einen Epimorphismus $I^{(\Lambda)}
\twoheadrightarrow M$ gibt. Faktormoduln und direkte Summen sind also wieder
$I$-generiert, ebenso -- wie aus der n"achsten Proposition folgt -- direkte
Produkte. Ein \emph{zyklischer} $R$-Modul $M$ ist genau dann $I$-generiert,
wenn es einen Epimorphismus $I \twoheadrightarrow M$ gibt.

\begin{Proposition}\label{1.1}
  F"ur einen $R$-Modul $M$ und ein Ideal $I$ von $R$ sind "aquivalent:
  \begin{iii}
    \item $M$ ist $I$-generiert.
    \item Es gibt eine Erweiterung $M \subset X$ mit $M = IX$.
  \end{iii}
   War $M \subset Y$ und $Y$ injektiv, so ist das weiter "aquivalent mit   
  \begin{iii}\setcounter{iiicount}{2}  
    \item $M = I(M :_{Y} I)$.
  \end{iii}
\end{Proposition}

\begin{Beweis}
(i $\to$ ii) Zum Epimorphismus $f \colon I^{(\Lambda)} \twoheadrightarrow M$
bilde man die Fasersumme
  \begin{equation*}
    \xymatrix{
      I^{(\Lambda)} \ar@{}[r]|{\textstyle\subset} \ar[d]_{f} & R^{(\Lambda)}
      \ar[d]^{g} \\
      M \ar@{}[r]|{\textstyle\subset} & X
    }
  \end{equation*}
und weil dann auch $g$ ein Epimorphismus ist, folgt $M = IX$.\\
(ii $\to$ i) Mit einem Epimorphismus $h \colon R^{(\Delta)}
\twoheadrightarrow X$ folgt $h(I \cdot R^{(\Delta)}) = IX = M$, also
$I^{(\Delta)} \twoheadrightarrow M$.\\
(iii $\to$ ii) ist klar, und bei (ii $\to$ iii) folgt aus dem kommutativen
Dreieck
\begin{equation*}
   \xymatrix{
      M \ar@{}[r]|{\textstyle\subset} \ar@{}[d]|{\textstyle\cap}  & X
      \ar@{-->}[dl]^{\alpha} \\
      Y & 
    }
\end{equation*}
mit $Y' = \operatorname{Bi} \alpha$, dass $M = I \cdot Y'$ ist, wegen $Y'
\subset M :_Y I$ also sogar $M = I(M :_Y I)$.
\end{Beweis}

Durch die Eigenschaft $M =IX$ ist $X$ nat"urlich nicht eindeutig bestimmt:
W"ahlt man eine reine Erweiterung $X \subset Y$ mit $I \cdot Y/X = 0$, folgt
auch $M =IY$. Andererseits kann man, wenn $M$ endlich erzeugt war, auch $X$
endlich erzeugt w"ahlen: $M = R a_1 + \ldots + R a_n$ und $a_i \in IX$, also
$a_i \in IX_i$ mit $X_i \subset X$ endlich erzeugt ($1 \leq i \leq n$)
$\mbox{}\implies X' = X_1 + \ldots + X_n$ endlich erzeugt und $M = IX'$. Mit
der minimalen Erzeugendenanzahl $v(M) = \dim_k(M/\mathfrak{m}
M)$  gilt genauer:

\begin{Lemma}\label{1.2}
  Ist $M$ endlich erzeugt und $I$-generiert, so gibt es eine Erweiterung $M
  \subset X$ mit $M = IX$ und $v(X) \leq v(M)$.
\end{Lemma}

\begin{Beweis}
  Sei $v(M) = n \geq 1$. Nach der Vorbemerkung kann man $X$ schon als
  endlich erzeugt annehmen, $X = X_1 + \ldots + X_t$, alle $X_i$ zyklisch
  $\neq 0$. Im $k$-Vektorraum $\overline{M} = M/\mathfrak{m} M =
  \overline{IX_1}+\ldots + \overline{IX_t}$ kann man so umordnen, dass die
  Darstellung $\overline{M} = \overline{IX_1}+\ldots + \overline{IX_s}$
  unverk"urzbar und $s \leq t$ ist. Es folgt $s \leq n$, $M = IX'$ mit $X' =
  X_1 + \ldots + X_s$, also $v(X') \leq s$.
\end{Beweis}

\noindent
\textbf{Beispiel 1} Genau dann ist $k$ $I$-generiert, wenn $I \neq 0$ ist.
Genau dann ist $R$ $I$-generiert, wenn $I \cong R$ ist.

\begin{Beweis}
  Das erste ist klar, und beim zweiten folgt aus $f \colon I
  \twoheadrightarrow R$, dass $I = \operatorname{Ke} f \oplus I_1$ ist, also
  $I_1 \cong R$ einen NNT enthält, also $\operatorname{Ke} f = 0$ ist.
\end{Beweis}

\noindent
\textbf{Beispiel 2} Sei $M$ ein beliebiger $R$-Modul, aber das Ideal $I$
zyklisch. Genau dann ist $M$ $I$-generiert, wenn $\operatorname{Ann}_R(I)
\cdot M = 0$ ist.

\begin{Beweis}
  "`$\Rightarrow$"' gilt stets nach (\ref{1.1}, ii), und bei
  "`$\Leftarrow$"' folgt aus $R^{(\Lambda)} \twoheadrightarrow M$ nach
  Voraussetzung auch $(R / \operatorname{Ann}_R(I))^{(\Lambda)}
  \twoheadrightarrow M$. Aber weil $I$ zyklisch war, ist $R /
  \operatorname{Ann}_R(I) \cong I$.
\end{Beweis}

\noindent
\textbf{Beispiel 3} Sei $R / \mathfrak{p}$ ein diskreter Bewertungsring, $I
\not\subset \mathfrak{p}$ und $\mathfrak{p} M = 0$. Dann ist $M$
$I$-generiert.

\begin{Beweis}
  Im Ring $\overline{R} = R / \mathfrak{p}$ ist das Ideal $\overline{I} = (I
  + \mathfrak{p}) / \mathfrak{p}$ nicht Null, also isomorph zu
  $\overline{R}$. Aus $\overline{I}^{(\Lambda)} \twoheadrightarrow M$ "uber
  $\overline{R}$ folgt dann $I^{(\Lambda)} \twoheadrightarrow M$ "uber $R$.
\end{Beweis}

\noindent
\textbf{Beispiel 4}
Sei $I$ von einer $R$-regul"aren Folge $(r_1, \dotsc, r_s)$ erzeugt ($s \geq
1$) und $IM = 0$. Dann ist $M$ $I$-generiert.

\begin{Beweis}
  Der Konormalenmodul $I/I^2$ ist "uber $R/I$ frei $\neq 0$ (siehe
  \cite{003} p.~153), so dass es einen Epimorphismus $I \twoheadrightarrow
  R/I$ gibt. Mit $R^{(\Lambda)} \twoheadrightarrow M$ folgt dann auch
  $I^{(\Lambda)} \twoheadrightarrow M$.
\end{Beweis}

\noindent
\textbf{Beispiel 5} Sei $M$ schwach-injektiv und $\operatorname{Ann}_R(I)
\cdot M = 0$. Dann ist $M$ $I$-generiert.

\begin{Beweis}
  $M$ ist in jeder Erweiterung koabgeschlossen und deshalb die kanonische
  Abb. $M/IM \longrightarrow M/M[\operatorname{Ann}_R(I)]$ ein wesentlicher
  Epimorphismus (\cite{006} p.~4402). Nach Voraussetzung ist $M =
  M[\operatorname{Ann}_R(I)]$, also sogar $M =IM$.
\end{Beweis}

Wir wollen im n"achsten Satz alle einreihigen, $I$-generierten $R$-Moduln
bestimmen. Dabei hei"se $M$ \emph{einreihig}, wenn sein Untermodulverband
$\mathcal{L}(M)$ totalgeordnet und endlich ist.

\begin{Lemma}\label{1.3}
  Ist $M$ von endlicher L"ange $n \geq 1$, so gilt: Genau dann ist $M$
  einreihig, wenn $\mathfrak{m}^{n-1} \not\subset \operatorname{Ann}_R(M)$.
\end{Lemma}

\begin{Beweis}
  Sei gleich $n \geq 2$. Bei "`$\Rightarrow$"' sei $M \supsetneqq U_1
  \supsetneqq U_2 \supsetneqq \ldots \supsetneqq U_{n-1} \supsetneqq 0$ der
  gesamte Untermodulverband. Dann sind alle $U_i$ zyklisch, also $U_1 =
  \mathfrak{m} M$, $U_2 = \mathfrak{m}U_1 = \mathfrak{m}^2M, \dotsc, U_{n-1}
  = \mathfrak{m}^{n-1}M \neq 0$.\\
  Bei "`$\Leftarrow$"' sind in der Folge $M \supset \mathfrak{m}M \supset
  \mathfrak{m}^2M \supset \ldots \supset \mathfrak{m}^{n-1}M \supsetneqq 0$
  sogar alle $\mathfrak{m}^iM \neq 0$ ($0 \leq i \leq n-1$), so dass in
  \begin{equation*}
    n = \text{L"ange}(M/\mathfrak{m}M) +
    \text{L"ange}(\mathfrak{m}M/\mathfrak{m}^2M) + \ldots +
    \text{L"ange}(\mathfrak{m}^{n-1}M) 
  \end{equation*}
  alle Summanden den Wert~1 haben, also alle $\mathfrak{m}^iM$ zyklisch
  sind. In $M \cong R/ \mathfrak{a}$ ist dann auch $\mathfrak{m} /
  \mathfrak{a}$ zyklisch, d.\,h. der Ring $R/\mathfrak{a}$ einreihig, also
  auch der $R$-Modul $M$.
\end{Beweis}

\begin{Satz}\label{1.4}
  F"ur einen einreihigen Modul $M$ der L"ange $n \geq 1$ sind "aquivalent:
  \begin{iii}
  \item $M$ ist $I$-generiert.
  \item $\mathfrak{a} = (\mathfrak{a} I) : I$ mit $\mathfrak{a} =
    \operatorname{Ann}_R(M)$.
  \item $\mathfrak{m}^{n-1} \cdot I \not\subset \operatorname{Ann}_R(M)
    \cdot I$.
  \end{iii}
\end{Satz}

\begin{Beweis}
  (i $\to$ ii) Es gen"ugt, dass $M$ zyklisch ist: Aus $I \twoheadrightarrow
  M \cong R/\mathfrak{a}$ folgt dann $I/\mathfrak{a}I \twoheadrightarrow
  R/\mathfrak{a}$, also $(\mathfrak{a}I) : I \subset \mathfrak{a}$, und
  "`$\supset$"' ist klar.\\
  (ii $\to$ i) Es gen"ugt, dass $R/\mathfrak{a}$ ein QF-Ring ist (d.\,h. als
  Ring artinsch und injektiv): Der nach Voraussetzung treue
  $R/\mathfrak{a}$-Modul $I/\mathfrak{a} I$ hat dann einen direkten
  Summanden $\cong R/\mathfrak{a}$, so dass aus $I \twoheadrightarrow
  I/\mathfrak{a} I \twoheadrightarrow R/\mathfrak{a}$ die Behauptung
  folgt.\\
  (ii $\to$ iii) W"are $\mathfrak{m}^{n-1}I \subset \mathfrak{a} I$, d.\,h.
  $\mathfrak{m}^{n-1} \subset (\mathfrak{a} I) : I$, folgte aus der
  Voraussetzung $\mathfrak{m}^{n-1} \subset \mathfrak{a}$ im Widerspruch zu
  (\ref{1.3}).\\
  (iii $\to$ ii) W"are $\mathfrak{a} \subsetneqq (\mathfrak{a} I) : I$, also
  $\mathfrak{m}^{n-1} + \mathfrak{a} \subset (\mathfrak{a} I) : I$, folgte
  mit $\mathfrak{m}^{n-1} I \subset \mathfrak{a} I$ ein Widerspruch zur
  Voraussetzung.
\end{Beweis}

\begin{Folgerung}\label{1.5}
  Sei $M$ einreihig von der L"ange $n \geq 1$ und $M \supsetneqq M_1
  \supsetneqq \ldots \supsetneqq M_n = 0$ seine Kompositionsreihe. Genau
  dann ist $M/M_i$ $I$-generiert ($1 \leq i \leq n$), wenn
  $\mathfrak{m}^{i-1}\cdot I \not\subset \operatorname{Ann}_R(M) \cdot I$.
  Insbesondere wird $M/M_i$ von $\mathfrak{m}^{n-i}$ generiert.
\end{Folgerung}

\begin{Beweis}
  Mit $M \cong R/\mathfrak{a}$ ist $M/M_i \cong R/\mathfrak{m}^i +
  \mathfrak{a}$ einreihig von der L"ange $i$, also nach (\ref{1.4}) genau
  dann $I$-generiert, wenn $\mathfrak{m}^{i-1} \cdot I \not\subset
  (\mathfrak{m}^i + \mathfrak{a}) \cdot I$ ist, mit $\overline{I} = I
  /\mathfrak{a} I$ also $\mathfrak{m}^{i-1} \cdot \overline{I} \not\subset
  \mathfrak{m}^i \cdot \overline{I}$, $\mathfrak{m}^{i-1} \cdot \overline{I}
  \neq 0$, $\mathfrak{m}^{i-1} \cdot I \not\subset \mathfrak{a}I$. Der
  Zusatz gilt, weil $\mathfrak{m}^{i-1} \cdot \mathfrak{m}^{n-i}$ nach
  (\ref{1.3}) nicht einmal in $\operatorname{Ann}_R(M)$ liegt.
\end{Beweis}

\section{Der Spezialfall $I = \mathfrak{m}$}

Bekanntlich ist $\mathcal{L}(R)$ genau dann totalgeordnet, wenn entweder $R$
ein diskreter Bewertungsring ist oder $R$ einreihig. Eine "ahnliche
Dichotomie erh"alt man f"ur die zyklischen $\mathfrak{m}$-generierten
$R$-Moduln:

\newcommand{\so}{\operatorname{So}}
\begin{Satz}\label{2.1}
  Ein zyklischer $R$-Modul $M$ ist genau dann $\mathfrak{m}$-generiert, wenn
  entweder $M \cong R / \mathfrak{p}$ ein diskreter Bewertungsring ist oder
  $M \cong C/\so(C)$ f"ur einen einreihigen $R$-Modul $C$. 
\end{Satz}

\begin{Beweis}
  "`$\Rightarrow$"' Sei $M$ zyklisch und $\mathfrak{m}$-generiert, $M =
  \mathfrak{m}X$, und nach (\ref{1.2}) kann man auch $X$ zyklisch annehmen.
  Zeigen wir zuerst, dass $\mathcal{L}(X)$ totalgeordnet ist: Bei $X/M = 0$,
  d.\,h. $M=0$ ist nichts zu zeigen. Bei $X/M \neq 0$ ist $X \cong
  R/\mathfrak{b}$, $M \cong \mathfrak{m}/\mathfrak{b}$, d.\,h. im Ring
  $R/\mathfrak{b}$ ist das maximale Ideal zyklisch. Damit ist der
  Idealverband von $R/\mathfrak{b}$ totalgeordnet, also auch
  $\mathcal{L}(X)$.\\
  Aus $M \cong R/\mathfrak{a}$ folgt entweder $\mathfrak{a} = \mathfrak{p}$
  und $R/\mathfrak{p}$ ein diskreter Bewertungsring. Oder $R/\mathfrak{a}$
  ist einreihig, aus $\mathfrak{m} \twoheadrightarrow R/\mathfrak{a}$,
  $\mathfrak{m}/\mathfrak{c} \cong R/\mathfrak{a}$ folgt dann, dass auch $C
  = R/\mathfrak{c}$ einreihig ist und $C/\so(C) \cong R/\mathfrak{c} :
  \mathfrak{m} = R /\mathfrak{a} \cong M$.\\
  "`$\Leftarrow$"' Falls $M \cong R/\mathfrak{p}$ ein diskreter
  Bewertungsring ist, folgt nach (1, Beispiel~3) die Behauptung. Falls $M
  \cong C/\so(C)$ ist mit $C$ einreihig, ist bei $M = 0$ wieder nichts zu
  zeigen, und bei $M \neq 0$ ist $\text{L"ange}(C) \geq 2$, also $C/\so(C)$
  nach (\ref{1.5}) $\mathfrak{m}$-generiert.
\end{Beweis}

\begin{Folgerung}
  Sei $M$ zyklisch und $\operatorname{Ann}_R(M) \subset \mathfrak{m}^2$.
  Genau dann ist $M$ $\mathfrak{m}$-generiert, wenn $\mathcal{L}(R)$
  totalgeordnet und $\so(R) \cdot M = 0$ ist.
\end{Folgerung}

\begin{Beweis}
  "`$\Rightarrow$"' Weil $\mathcal{L}(M)$ totalgeordnet ist, also auch
  $\mathcal{L}(\mathfrak{m}/\mathfrak{m}^2)$, ist $\mathfrak{m}$ zyklisch,
  also sogar $\mathcal{L}(R)$ totalgeordnet, und $\so(R) \cdot M = 0$ ist
  klar.\\
  "`$\Leftarrow$"' Dann ist entweder $R$ ein diskreter Bewertungsring und
  nach (1, Beispiel~3) jeder $R$-Modul $\mathfrak{m}$-generiert, oder $R$
  einreihig und $M$ nach (1, Beispiel~2) $\mathfrak{m}$-generiert.
\end{Beweis}

In dem Spezialfall, dass $(R, \mathfrak{m})$ ein QF-Ring, aber nicht
einreihig ist, $\mathfrak{m}^2 \neq 0$ und $\mathfrak{m}^3 = 0$ (z.\,B. $R =
k[X,Y]/(X^2,Y^2)$), wollen wir jetzt alle $\mathfrak{m}$-generierten Unter-
und Faktormoduln von $R$ bestimmen:

\medskip

\noindent
\textbf{Beispiel}
  Ist $R$ wie eben, so sind $0$, $\mathfrak{m}^2$, $\mathfrak{m}$ die
  einzigen $\mathfrak{m}$-generierten Untermoduln von $R$ und $R/R$,
  $R/\mathfrak{m}$ die einzigen $\mathfrak{m}$-generierten Faktormoduln von
  $R$.

\begin{Beweis}
  Klar sind die angegebenen Moduln alle $\mathfrak{m}$-generiert.
  Untermoduln: Ist $\mathfrak{a} \subsetneqq \mathfrak{m}$
  $\mathfrak{m}$-generiert, folgt (weil ${}_RR$ injektiv ist) nach
  (\ref{1.1}, iii) $\mathfrak{a} = \mathfrak{m} \mathfrak{b}$ f"ur ein Ideal
  $\mathfrak{b} \neq R$, also $\mathfrak{a} \in \{0, \mathfrak{m}^2\}$.\\
  Faktormoduln: $R$ besitzt keinen einreihigen Modul $M$ der L"ange~3 (also
  auch keinen gr"o"serer L"ange), denn sonst folgte nach (\ref{1.3})
  $\mathfrak{m}^2 \not\subset \operatorname{Ann}_R(M)$, also $M \cong R$
  entgegen der Voraussetzung. Nach (\ref{2.1}) kann deshalb
  $R/\mathfrak{a}$, mit $\mathfrak{a} \subsetneqq \mathfrak{m}$ nicht
  $\mathfrak{m}$-generiert sein.
\end{Beweis}

\section{Dualisierung}

Verlangt man von einem $R$-Modul $M$, dass nicht er selbst, sondern sein
Matlis-Duales $M^{\circ} = \operatorname{Hom}_R(M,E)$ $I$-generiert ist,
erh"alt man eine neue Klasse von $R$-Moduln, die im Spezialfall $I =
\mathfrak{m}$ zu den von Heinzer, Ratliff~Jr. und Rush in \cite{002}
untersuchten "`basically full ideals"' f"uhrt.

\begin{Proposition}
  F"ur einen $R$-Modul $M$ und ein Ideal $I$ von $R$ sind "aquivalent:
  \begin{iii}
  \item $M^{\circ}$ ist $I$-generiert.
  \item $M$ ist $I^{\circ}$-kogeneriert.
  \item Es gibt einen $R$-Modul $C$ mit $C/C[I] \cong M$.
  \end{iii}
  War $M \cong A/B$ und $A$ frei, so ist das weiter "aquivalent mit
  \begin{iii}\setcounter{iiicount}{3}
  \item $B = (IB) :_A I$.
  \end{iii}
\end{Proposition}

\begin{Beweis}
  (i $\to$ ii) Aus $I^{(\Lambda)} \twoheadrightarrow M^{\circ}$ folgt $M
  \hookrightarrow M^{\circ\circ} \hookrightarrow (I^{(\Lambda)})^{\circ}
  \cong (I^{\circ})^{\Lambda}$.\\
  (ii $\to$ i) F"ur jeden $I$-generierten $R$-Modul $N = I \cdot Y$ ist
  $N^{\circ} \cong Y^{\circ}/Y^{\circ}[I]$, also auch $N^{\circ\circ} \cong
  I \cdot Y^{\circ\circ}$ $I$-generiert. Unter unserer Voraussetzung, d.\,h.
  $M \hookrightarrow (I^{\circ})^{\Lambda} \cong N^{\circ}$ mit $N =
  I^{(\Lambda)}$, folgt jetzt $N^{\circ\circ} \twoheadrightarrow M^{\circ}$,
  und weil $I$, also auch $N$ und $N^{\circ\circ}$ $I$-generiert sind, ist
  es auch $M^{\circ}$.\\
  (i $\to$ iii) Aus $M^{\circ} = I \cdot Y$ folgt $\beta\colon Y^{\circ}
  \twoheadrightarrow M^{\circ\circ}$ mit $\operatorname{Ke} \beta =
  \operatorname{Ann_{Y^{\circ}}}(IY) \cong Y^{\circ}[I]$, mit $C =
  \beta^{-1}(M)$ und der induzierten Abb. $\gamma\colon C \twoheadrightarrow M$
  also $\operatorname{Ke} \gamma = C[I]$.\\
  (iii $\to$ i) Es folgt $M^{\circ} \cong \operatorname{Ke}(C^{\circ}
  \twoheadrightarrow (C[I])^{\circ}) = \operatorname{Ann}_{C^{\circ}}(C[I])
  = I \cdot C^{\circ}$.\\
  (iii $\to$ iv) Weil $A$ frei ist, erh"alt man ein kommutatives Diagramm
  \begin{equation*}
    \xymatrix{
      & & & A \ar@{-->}[dl]_{\alpha} \ar[d]^{\mathrm{kan}}& \\
      0 \ar[r] & C[I] \ar@{}[r]|{\textstyle\subset} & C \ar[r]_{\gamma} &
      A/B \ar[r] & 0\ ,
    }
  \end{equation*}
  in dem $B = \alpha^{-1}(\operatorname{Ke} \gamma) =
  \operatorname{Ke}\alpha :_A I$ ist, also $B = (IB) :_A I$.\\
  (iv $\to$ iii) Hier kann $A$ beliebig sein, denn mit $C = A/IB$ gilt
  $C/C[I] \cong A/(IB) :_A I = A/B \cong M$.
\end{Beweis}

Auch f"ur einen zyklischen, $I$-generierten $R$-Modul $M \cong
R/\mathfrak{b}$ gilt nach dem Beweis von (\ref{1.4}, i $\to$ ii), dass
$\mathfrak{b} = (\mathfrak{b}I) : I$ ist.

\begin{Folgerung}\label{3.2}
  Jeder zyklische $I$-generierte $R$-Modul ist auch $I^{\circ}$-kogeneriert.
\end{Folgerung}

Ein $R$-Modul $M$ hei"st bekanntlich \emph{kozyklisch}, wenn er sich in $E$
einbetten l"asst. Insbesondere ist dann $M^{\circ} / \mathfrak{m} \cdot M^{\circ}
\cong (M[\mathfrak{m}])^{\circ}$ Null oder einfach. War also $M$ sowohl
zyklisch als auch kozyklisch, ist $M$ von endlicher L"ange, also auch
$M^{\circ}$ zyklisch und daher $M \cong M^{\circ}$.

\begin{Folgerung}\label{3.3}
  Ist $M$ sowohl zyklisch als auch kozyklisch, so gilt: Genau dann ist $M$
  $I^{\circ}$-kogeneriert, wenn $M$ $I$-generiert ist.
\end{Folgerung}

\noindent
\textbf{Beispiele} f"ur $I^{\circ}$-kogenerierte $R$-Moduln: (1) Eine
direkte Summe von einreihigen Moduln ist genau dann $I^{\circ}$-kogeneriert,
wenn sie $I$-generiert ist (denn die unzerlegbaren Bausteine erf"ullen
(\ref{3.3})). (2) Ist $R$ ein diskreter Bewertungsring und $I \neq 0$, so
ist jeder $R$-Modul $M$ $I^{\circ}$-kogeneriert (denn $M^{\circ}$ ist
$I$-generiert). (3) Wird $I$ von einer $R$-regul"aren Folge $(r_1, \dotsc,
r_s)$ erzeugt ($s \geq 1$) und ist $IM = 0$, so ist $M$
$I^{\circ}$-kogeneriert (denn es ist $I \cdot M^{\circ} = 0$). (4) $E$ ist
genau dann $I^{\circ}$-kogeneriert, wenn $I \cong R$ ist (denn aus $E
\hookrightarrow I^{\circ} \cong E/E[I]$ folgt $E \cong I^{\circ}$, $I \cong
R$). (5) Ist $R$ ein QF-Ring wie im 2.~Abschnitt, d.\,h. nicht einreihig
und $\mathfrak{m}^2 \neq 0$, $\mathfrak{m}^3 = 0$, so sind $0$,
$\mathfrak{m}^2$ die einzigen $\mathfrak{m}^{\circ}$-kogenerierten Untermoduln
und $R/R$, $R/\mathfrak{m}$, $R/\mathfrak{m}^2$ die einzigen
$\mathfrak{m}^{\circ}$-kogenerierten Faktormoduln von $R$.

\bigskip

Ist $A$ ein beliebiger $R$-Modul, so hei"se ein Untermodul $B$
\emph{$I$-gro"s} in $A$, wenn es ein $B' \subset B$ gibt mit $B/B' =
(A/B')[I]$. Das ist "aquivalent mit $B = (IB) :_A I$, und im Spezialfall $I
= \mathfrak{m}$, $A$ endlich erzeugt, $A/B$ artinsch und $B \neq 0$ ist
das gerade eine Beschreibung von "`$B$ basically full in $A$"' in (\cite{002}
Theorem~2.12). Sind aber $B \subset A$ und $I$ beliebig, gilt immer noch:

\begin{myenum}
\item $B^{*} = (IB) :_A I$ ist der kleinste Modul zwischen $B$ und $A$, so
  dass $B^{*}$ $I$-gro"s in $A$ ist (\cite{002} Theorem~4.2).
\item Sei $B$ $I$-gro"s in $A$ und $(A/B)[I]$ gro"s in $A/B$. Dann folgt aus
  $B \subset C \subset A$ und $B \cap IC = IB$ stets $B=C$ (\cite{002}
  Theorem~2.12).
\item Aus $B \subset C \subset A$ und $B \cap IC = IB$ folge stets $B = C$.
  Falls dann $A$ endlich erzeugt ist, wird $A/B$ durch eine Potenz von $I$
  annulliert (\cite{002} Theorem~2.6).
\end{myenum}

Die Beweise sind direkte Verallgemeinerungen von denen in \cite{002}, sollen
aber im dualen Fall im n"achsten Abschnitt ausgef"uhrt werden.

\section{$I$-kleine Erweiterungen}

Sei $M \subset X$ eine beliebige Erweiterung. Wir sagen, $M$ sei
\emph{$I$-klein} in $X$, wenn es einen Zwischenmodul $M \subset X' \subset
X$ gibt mit $M = IX'$, und daraus folgt sofort $M = I(M :_X I)$. Nat"urlich
ist dann $M$ $I$-generiert, und falls $X$ injektiv ist, gilt nach
(\ref{1.1}, iii) auch die Umkehrung.

\begin{Lemma}\label{4.1}
  Sei $M \subset X$ eine beliebige Erweiterung und $M_{*} = I(M : I)$.
  \begin{abc}
  \item Stets gilt $IM \subset M_{*} \subset M$ und $M_{*} : I = M : I$.
  \item $M_{*}$ ist der gr"o"ste Untermodul von $M$, der $I$-klein in $X$ ist.
  \end{abc}
\end{Lemma}

\begin{Beweis}
  F"ur Untermoduln $W$ von $X$ schreiben wir statt $W :_X I$ nur mehr $W :
  I$.\\
  (a) Aus $M \subset M : I$ folgt $IM \subset M_{*}$. Aus $u \in M_{*}$
  folgt $u = \sum r_ix_i$ mit $r_i \in I$, $x_i \in M : I$, also $u \in M$.
  Aus $x \in M : I$ folgt $rx \in M_{*}$ f"ur alle $r \in I$, also $x \in
  M_{*} : I$.\\
  (b) Multipliziert man die letzte Gleichung mit $I$, erh"alt man $M_{**} =
  M_{*}$, und das bedeutet, dass $M_{*}$ $I$-klein in $X$ ist. Falls auch $U
  \subset M$ und $U$ $I$-klein in $X$ ist, folgt aus $U_{*} = U$ sogar $U
  \subset M_{*}$.
\end{Beweis}

\begin{Bemerkung}\label{4.2}
  (1) F"ur jede Erweiterung $M \subset X$ gilt: $M =IX \Rightarrow M$ ist
  $I$-klein in $X \Rightarrow M \subset IX$. Das nachfolgende Beispiel (mit
  $X = R$ und $I = \mathfrak{m}$) zeigt, dass nirgends die Umkehrung gilt.
  (2) Ist $M \subset X$ beliebig und $I$ zyklisch, gilt $M_{*} = M \cap IX$.
  Insbesondere ist $M$ genau dann $I$-klein in $X$, wenn $M \subset IX$ ist.
  (3) Ist $M \subset X$ eine reine Erweiterung, folgt aus $X[I]
  \twoheadrightarrow \frac{X}{M}[I]$, dass $M_{*} = IM$ ist. Insbesondere
  ist $M$ genau dann $I$-klein in $X$, wenn $M$ $I$-teilbar ist.
\end{Bemerkung}

\noindent
\textbf{Beispiel} Sei $R$ ein QF-Ring wie im 2.~Abschnitt, d.\,h. nicht
einreihig und $\mathfrak{m}^2 \neq 0$, $\mathfrak{m}^3 = 0$. F"ur jedes
Ideal $\mathfrak{b} \subset R$ l"asst sich dann $\mathfrak{b}_{*} =
\mathfrak{m}(\mathfrak{b} : \mathfrak{m})$ und $\mathfrak{b}^{*} =
(\mathfrak{m}\mathfrak{b}) : \mathfrak{m}$ berechnen:

\smallskip

\begin{equation*}
  \begin{array}{l|ccccc}
    \mathfrak{b} & 0 & \mathfrak{m}^2 & \ \mathfrak{m}^2 \subsetneqq
    \mathfrak{b} \subsetneqq \mathfrak{m}\ & \mathfrak{m} & R \\[1ex] \hline
    \rule{0ex}{3ex}\mathfrak{b}_{*} & 0 & \ \mathfrak{m}^2\   &
    \mathfrak{m}^2 & \ \mathfrak{m}\  &
    \ \mathfrak{m}\  \\[1ex]
    \mathfrak{b}^{*}\  & \ \mathfrak{m}^2\ & \ \mathfrak{m}^2\  & \mathfrak{m}^{\phantom{2}} &
    \mathfrak{m} & R \\[1ex] \hline
  \end{array}
\end{equation*}

\medskip

\begin{Lemma}\label{4.3}
  Genau dann ist $M$ $I$-klein in $X$, wenn
  $\operatorname{Ann}_{X^{\circ}}(M)$ $I$-gro"s in $X^{\circ}$ ist.
\end{Lemma}

\begin{Beweis}
  F"ur eine beliebige Erweiterung $C \subset A$ gilt
  $\operatorname{Ann}_{A^{\circ}}(I \cdot C) =
  \operatorname{Ann}_{A^{\circ}}(C) :_{A^{\circ}} I$ und
  $\operatorname{Ann}_{A^{\circ}}(C :_A I) = I \cdot
  \operatorname{Ann}_{A^{\circ}}(C)$. Speziell f"ur $A = X^{\circ}$ und $B =
  \operatorname{Ann}_{X^{\circ}}(M)$ gilt damit
  \begin{equation*}B^{*} = (IB) :_A I = (I \cdot
    \operatorname{Ann}_{X^{\circ}}(M)) :_{X^{\circ}} I =
    \operatorname{Ann}_{X^{\circ}}(M : I) :_{X^{\circ}} I =
    \operatorname{Ann}_{X^{\circ}}(M_{*}),
   \end{equation*} also $M_{*} = M$ genau dann,
    wenn $B = B^{*}$ ist.
\end{Beweis}

\begin{Bemerkung}\label{4.4}
  Beginnt man mit einer beliebigen Erweiterung $B \subset A$, so kann man
  entsprechend zeigen: Genau dann ist $B$ $I$-gro"s in $A$, wenn
  $\operatorname{Ann}_{A^{\circ}}(B)$ $I$-klein in $A^{\circ}$ ist.
\end{Bemerkung}

$M$ hei"se \emph{stark $I$-klein} in $X$, wenn f"ur jeden Untermodul $U$ von
$M$ gilt: Ist $\frac{X}{U}[I] \to \frac{X}{M}[I]$ surjektiv (d.\,h. $U : I +
M = M : I$), so folgt $U = M$. Speziell f"ur $U = M_{*}$ ist diese Gleichung
nach (\ref{4.1}, a) erf"ullt, also nach Voraussetzung $M_{*} = M$: Jede
stark $I$-kleine Erweiterung ist $I$-klein. Die Umkehrung gilt i.~Allg.
nicht: Ist $M$ stark $I$-klein in $X$ und entweder $M$ rein in $X$ oder
$(X/M)[I] = 0$, folgt bereits $M = 0$.

\begin{Lemma}\label{4.5}
  Ist $M$ stark $I$-klein in $X$ und $X_1 \subset X$, so ist auch $(M
  +X_1)/X_1$ stark $I$-klein in $X/X_1$.
\end{Lemma}

\begin{Beweis}
  Sei $\overline{X} = X/X_1$ und $\overline{V} \subset \overline{M}$, so
  dass $\frac{\overline{X}}{\overline{V}}[I] \twoheadrightarrow
  \frac{\overline{X}}{\overline{M}}[I]$ surjektiv ist, also $X_1 \subset V
  \subset M + X_1$ und $\frac{X}{V}[I] \twoheadrightarrow
  \frac{X}{M+X_1}[I]$ surjektiv, d.\,h.
  \begin{equation*}
    V : I + (M+X_1) = (M+X_1) : I\ .
  \end{equation*}
  K"onnten wir $(V \cap M) : I + M = M : I$ zeigen, folgte nach
  Voraussetzung $V \cap M = M$, $M \subset V$, also $\overline{V} =
  \overline{M}$ wie verlangt.\\
  F"ur $x \in M : I$ gilt $x = a + b$ mit $a \in V : I$, $b \in M$, also $ra
  = rx - rb \in V \cap M$ f"ur alle $r \in I$, also $a \in (V \cap M) : I$,
  $x \in (V \cap M) : I + M$ wie behauptet.  
\end{Beweis}

Bekanntlich ist $IM$ genau dann klein in $M$, wenn $I \subset \bigcap
\operatorname{Koass}(M)$ ist (\cite{004} Lemma~2.2). Das vererbt sich auf
Faktormoduln, und war $M/U$ sogar artinsch, folgt aus $I \subset
\sqrt{\operatorname{Ann}_R(M/U)}$, dass $I^n \cdot M/U  = 0$ ist f"ur ein
$n \geq 1$. Umgekehrt ist, falls jeder artinsche Faktormodul von $M$ durch
eine Potenz von $I$ annulliert wird, auch $IM$ klein in $M$: Aus $V + IM=M$
mit $M/V$ artinsch folgt $I^n \cdot M/V = 0$, $V + I^nM = M$, $V = M$. 

\begin{Satz}\label{4.6}
  Genau dann ist $M$ stark $I$-klein in $X$, wenn $M$ $I$-klein in $X$ ist
  und $IM$ klein in $M$.
\end{Satz}

\begin{Beweis}
  "`$\Leftarrow$"' ist klar: Aus $U \subset M$ und $U : I + M = M : I$ folgt
  durch Multiplikation mit $I$, dass $U_{*} + IM = M_{*}$ ist, mit der
  ersten Bedingung also $U_{*} + IM = M$, mit der zweiten dann $U = M$.\\
  Bei "`$\Rightarrow$"' nehmen wir im \textbf{1.~Schritt} an, dass $X$
  zus"atzlich artinsch ist. Nach der dualen Formulierung von Artin-Rees
  (\cite{005} p.~640) gibt es dann ein $n \geq 1$ mit
  \begin{equation*}
    \frac{X}{M}[I]\; \subset \; \frac{X[I^n]+M}{M}\ .
  \end{equation*}
  Mit $U = M[I^n]$ folgt $U : I + M = M : I$ (also nach Voraussetzung $U =
  M$, $I^n \cdot M = 0$ wie gewünscht): F"ur $x \in M : I$ gilt $x = a+ b$
  mit $a \in X[I^n]$, $b \in M$, also $ra = rx-rb \in M \cap X[I^n] = U$
  f"ur alle $r \in I$, also $a \in U : I$, $x \in U : I + M$ wie
  behauptet.\\
  Im \textbf{2.~Schritt} sei jetzt $X$ beliebig. Nach der Vorbemerkung
  m"ussen wir zeigen, dass jeder artinsche Faktormodul $M/V$ durch eine
  Potenz von $I$ annulliert wird. Ist $X_1/V$ ein Durchschnittskomplement
  von $M/V$ in $X/V$, wird $M/V \hookrightarrow X/X_1$ ein wesentlicher
  Monomorphismus, so dass auch $\overline{X} = X/X_1$ artinsch und
  $\overline{M}$ nach (\ref{4.5}) stark $I$-klein in $\overline{X}$ ist.
  Nach dem ersten Schritt folgt jetzt $I^n \cdot \overline{M} = 0$ f"ur ein
  $n \geq 1$, d.\,h. $I^n \cdot M/V = 0$.
\end{Beweis}

\begin{Bemerkung}
  Ein Untermodul $B$ von $A$ hei"se \emph{stark $I$-gro"s} in $A$, wenn aus
  $B \subset C \subset A$ und $B \cap IC = IB$ stets folgt $B=C$. Die zu
  (\ref{4.6}) duale Aussage lautet dann: Genau dann ist $B$ stark $I$-gro"s
  in $A$, wenn $B$ $I$-gro"s in $A$ ist und $\frac{A}{B}[I]$ gro"s in
  $\frac{A}{B}$ (d.\,h. $A/B$ $I$-torsion ist). Das ist eine Zusammenfassung der
  Punkte~(2) und (3) im dritten Abschnitt.
\end{Bemerkung}

\begin{Folgerung}
  Genau dann ist $M$ stark $I$-klein in $X$, wenn
  $\operatorname{Ann}_{X^{\circ}}(M)$ stark $I$-gro"s in $X^{\circ}$ ist.
\end{Folgerung}

\begin{Beweis}
  Auf Grund von (\ref{4.3}) ist nur noch zu zeigen: Genau dann ist $IM$
  klein in $M$, wenn $X^{\circ}/\operatorname{Ann}_{X^{\circ}}(M) \cong
  M^{\circ}$ $I$-torsion ist. Das ist aber klar wegen $\operatorname{Koass}(M) =\operatorname{Ass}(M^{\circ})$.
\end{Beweis}

\begin{Bemerkung}
  Entsprechend (\ref{4.4}) kann man zeigen: Ist $B \subset A$ und
  $\operatorname{Ann}_{A^{\circ}}(B)$ stark $I$-klein in $A^{\circ}$, so ist
  $B$ stark $I$-gro"s in $A$. Die Umkehrung gilt aber nicht: Ist $R$ ein
  diskreter Bewertungsring und $A = K^{(\mathbb{N})}$, $B =
  R^{(\mathbb{N})}$, so ist $B$ stark $\mathfrak{m}$-gro"s in $A$, aber
  $\operatorname{Ann}_{A^{\circ}}(B)$ \emph{nicht} stark
  $\mathfrak{m}$-klein in $A^{\circ}$, denn in
  $\operatorname{Ann}_{A^{\circ}}(B) \cong \widehat{R}^{\mathbb{N}}$ ist das
  Radikal nicht klein.
\end{Bemerkung}

\end{document}